\documentclass[12pt, reqno]{amsart}
\usepackage{amsmath,amsopn,amssymb,amsthm}
\usepackage[english]{babel}
\usepackage[backref=page, breaklinks=true,colorlinks=true,linkcolor=blue,citecolor=blue,urlcolor=blue]{hyperref}
\usepackage{graphicx}

\renewenvironment{proof}[1][Proof]{\textbf{#1.} }{\ \rule{0.5em}{0.5em}}

\DeclareMathOperator{\length}{length}
\DeclareMathOperator{\per}{per}
\DeclareMathOperator{\area}{area}
\DeclareMathOperator{\vol}{vol}
\DeclareMathOperator{\diam}{diam}
\DeclareMathOperator{\bd}{bd}
\DeclareMathOperator{\co}{co}

\DeclareMathOperator{\intt}{int}

\newtheorem{theorem}{Theorem}

\newtheorem{corollary}{Corollary}

\theoremstyle{definition}

\newtheorem{remark}{Remark}

\newtheorem{problem}{Problem}

\begin{document}

\title
[Properties of a curve whose convex hull \dots]
{Properties of a curve whose convex hull \\covers a given convex body}
\author{Yu.G.~Nikonorov}

\address{Nikonorov\ Yuri\u\i\  Gennadievich\newline
Southern Mathematical Institute of \newline
the Vladikavkaz Scientific Center of \newline
the Russian Academy of Sciences, \newline
Vladikavkaz, Markus st., 22, \newline
362027, RUSSIA}
\email{nikonorov2006@mail.ru}

\begin{abstract}
In this note, we prove the following inequality for the norm of a convex body $K$ in $\mathbb{R}^n$, $n\geq 2$:
$N(K) \leq \frac{\pi^{\frac{n-1}{2}}}{2 \Gamma \left(\frac{n+1}{2}\right)}\cdot \length (\gamma) +
\frac{\pi^{\frac{n}{2}-1}}{\Gamma \left(\frac{n}{2}\right)} \cdot  \diam(K)$,
where $\diam(K)$ is the diameter of $K$,
$\gamma$ is any curve in $\mathbb{R}^n$ whose convex hull covers~$K$, and $\Gamma$ is the gamma function.
If  in addition $K$ has constant width $\Theta$,
then we get the inequality $\length (\gamma) \geq \frac{2(\pi-1)\Gamma \left(\frac{n+1}{2}\right)}{\sqrt{\pi}\,\Gamma \left(\frac{n}{2}\right)}\cdot \Theta
\geq 2(\pi-1) \cdot \sqrt{\frac{n-1}{2\pi}}\cdot \Theta$.
In addition, we pose several unsolved problems.

\vspace{2mm}

\vspace{2mm} \noindent Key words and phrases:  convex body, convex hull, curve, diameter, quermassintegral, norm of a convex body.
\end{abstract}

\maketitle

\section{Introduction and main results}\label{sect.1}

There are many unsolved problem related to the convex hull of a curve in Euclidean space, see e.~g. \cite{Ghomi2018} and \cite{Zal1996} for discussions of some of them.
In this note, we consider some relationships between a convex body $K$ in Euclidean space and a curve $\gamma$, whose convex hull cover $K$.
We refer to \cite{BBI2001} for standard results in the metric geometry and to \cite{BoFe1987, Had1957, MaMoOl} for classical results in the geometry of convex bodies.

We identify $n$-dimensional Euclidean space with $\mathbb{R}^n$ supplied with the standard Euclidean metric~$d$, where $d(x,y)=\sqrt{\sum\limits_{i=1}^n (x_i-y_i)^2}$.
For any subset $A\subset \mathbb{R}^n$, $\co (A)$ means the convex hull of $A$. For every points $B,C \in \mathbb{R}^n$, $[B,C]$
denotes the line segment between these points.

{\it A convex body} is any compact convex subset of $\mathbb{R}^n$ (convex bodies with empty interior are allowed).
We shall denote by $\vol(K)$, $\area(K)$, $\bd (K)$ and $\intt (K)$ respectively the volume, the surface area, the boundary, and the interior of a convex figure $K$.
Note also that the diameter $\diam (K):=\max \left\{ d(x,y)\,|\, x,y \in K \right\}$
of a convex body $K$ coincides with the maximal distance between two parallel support lines
to $K$.

The symbol $B^n(x,\rho)$ denotes a closed ball in $\mathbb{R}^n$ with center $x \in  \mathbb{R}^n$ and radius $\rho \geq 0$.
We consider also the unit ball $B^n:=B^n(0,1)$ and the unit sphere $S^{n-1}=\bd (B^n)$.
We use the symbols $\omega_n$ and $\sigma_{n-1}$ for the volume of the unit $B^n$ and for the surface area of the unit sphere $S^{n-1}$ respectively.
Recall that $\sigma_{n-1}=n \omega_n=\frac{2\pi^{n/2}}{\Gamma(n/2)}$, where $\Gamma$ is the gamma function.

{\it A curve} $\gamma$ is the image of a continuous mapping $\varphi:[a,b]\subset \mathbb{R} \mapsto \mathbb{R}^n$.
As usually, the length of $\gamma$ is defined as  $\length (\gamma):=\sup \left\{ \sum\limits_{i=1}^m d(\varphi(t_{i-1}),\varphi(t_i))\right\}$,
where the supremum is taken over all finite increasing sequences $a=i_0<i_1<\cdots <i_{m-1}<i_m=b$ that lie in
the interval $[a,b]$.
A curve $\gamma$ is called {\it rectifiable} if $\length (\gamma) <\infty$.

For $n=2$, we use the term {\it a convex figure} for a convex body $K\subset \mathbb{R}^2$, $\area(K)$ is the perimeter $\per(K)$ of $K$ in this case.
We call a curve $\gamma\subset \mathbb{R}^2$ {\it convex} if it is a closed connected subset of the boundary
of the convex hull $\co(\gamma)$ of $\gamma$.
\medskip

Let us recall a remarkable property of planar curves, that was obtained in \cite{NikNik2021}.

\begin{theorem}[\cite{NikNik2021}]\label{the1}
For a given convex figure $K$ and for any planar curve $\gamma$ with the property $K\subset \co(\gamma)$, the inequality
\begin{equation}\label{eq.main}
\length (\gamma)\geq \per (K) - \diam(K)
\end{equation}
holds. Moreover, this inequality becomes an equality if and only if
$\gamma$ is a convex curve, $\bd (K)=\gamma \cup [A,B]$, and $\diam(K)=d(A,B)$, where $A$ and $B$ are the endpoints of $\gamma$.
\end{theorem}

The inequality \eqref{eq.main} was suggested by A.~Akopyan and V.~Vysotsky, who proved it in the case when $\gamma$ is passing through
all extreme points of $K$ (see Theorem 7 in \cite{AkVys2017}).

\medskip

In this note, we obtain some generalization of Theorem \ref{the1} to the multidimensional case.
Since $\length(\gamma)$ and $\diam(K)$ are respectively in some sense ``one-dimensional'' characteristics of a curve $\gamma$ and a convex body $K$ in $\mathbb{R}^n$,
then we should find some ``one-dimensional'' characteristic of the convex body $K\subset \mathbb{R}^n$ instead of $\per (K)$, used in the two-dimensional case.
Such a natural generalization of $\per (K)$ is {\it the norm} $N(K)$ of convex body $K$.
\smallskip

Let $K\subset \mathbb{R}^n$  be a non-empty compact convex set. The support function
$h(K, \cdot)$ of $K$ is defined by
$$
h(K, u) := \sup \{(x,u)\,|\, x\in K \}\,.
$$ for $u\in \mathbb{R}^n$.
The width function $w(K,\cdot)$ of $K$
is defined by
$$
w(K,u):= h(K,u)+h(K,-u), \qquad  u \in  S^{n-1}\,.
$$
The mean value of the width function over $S^{n-1}$ is called {\it the mean width} and
denoted by $w(\cdot)$, thus
$w(K) :=\frac{1}{\sigma_{n-1}} \int_{S^{n-1}} h(K,u) du$. {\it The norm} of $K$ could be defined as
\begin{equation}\label{eq.norm}
N(K):=\frac{1}{2} \int_{S^{n-1}} h(K,u) du=\frac{\sigma_{n-1}}{2} \cdot w(K)\,.
\end{equation}
See details e.~g. in 6.1.7 of \cite{Had1957}.

For a convex body $K \subset \mathbb{R}^n$ and for a number $\varepsilon \geq 0$, we consider
{\it the outer parallel body of $K$ at distance $\varepsilon$}, which is defined as follows:
$$
K_{\varepsilon} := \{x \in \mathbb{R}^n\,|\, d(K,x)\leq \varepsilon\}=\bigcup\limits_{x\in K} B^n(x,\varepsilon).
$$

{\it The Steiner formula} says that its volume can be expressed
as a polynomial of degree $n$ in the parameter  $\varepsilon$:
\begin{equation}\label{eq.steiner}
\vol(K_{\varepsilon})=\sum\limits_{k=0}^n C_n^k\,\cdot W_k(K) \,\cdot \varepsilon^k, \qquad C_n^k=\frac{n!}{k! (n-k)!}.
\end{equation}
The functionals $W_0, W_1, \dots, W_{n-1}, W_n$ are called the quermassintegrals, see \cite{BoFe1987, Had1957}.
Note that $W_0(K)=\vol(K)$, $nW_1(k)=\area(K)$, $W_n=\omega_n$ and $nW_{n-1}=N(K)$.
In particular, we have $N(K)=\area(K)=\per(K)$ for $n=2$.

\medskip

The main result of this note is the following theorem, that is proved in the next section.

\begin{theorem}\label{the2}
For a given convex body $K \subset \mathbb{R}^n$ and for any  curve $\gamma \subset \mathbb{R}^n$ with the property $K\subset \co(\gamma)$, where $ n \geq 2$,
the inequality
\begin{equation}\label{eq.main.new1}
N(K)=nW_{n-1}(K) \leq \frac{\pi^{\frac{n-1}{2}}}{2 \Gamma \left(\frac{n+1}{2}\right)}\cdot \length (\gamma) +
\frac{\pi^{\frac{n}{2}-1}}{\Gamma \left(\frac{n}{2}\right)} \cdot  \diam(K)
\end{equation}
holds.
\end{theorem}

\begin{remark}\label{rem1}
The inequality \eqref{eq.main.new1} is equivalent to \eqref{eq.main} for $n=2$ and
does not hold in general for $n=1$. For $n=3$, \eqref{eq.main.new1} has the form $N(K)\leq \frac{\pi}{2} \length (\gamma) +
2 \diam(K)$.
\end{remark}

It should be noted that in the equality in \eqref{eq.main.new1} does not hold if $n \geq 3$ and $\diam(K) >0$ (i.~e. $K$ is not a one-point set),
see Remark \ref{rem3}.
\smallskip

By virtue of the assertion of Theorem \ref{the2}, the following problem seems to be interesting.

\begin{problem}\label{pr1}
For a given $n$, find all possible constant $A(n)$ and $B(n)$ such that
$$
N(K) \leq A(n)\cdot \length (\gamma) +
B(n) \cdot  \diam(K)
$$
for any convex body $K \subset \mathbb{R}^n$ and for any  curve $\gamma \subset \mathbb{R}^n$ with the property $K\subset \co(\gamma)$.
\end{problem}

\medskip

Now, let us suppose in addition that in conditions of Theorem \ref{the2} the convex body $K$ has constant width $\Theta$.
Then, by the definition of the mean width
$w(K)$
we have $w(K)=\Theta$ (since $w(K,u) \equiv \Theta$ for all $u\in S^{n-1}$), and by the definition of norm $N(K)$ (see \eqref{eq.norm}) we have
$N(K)=\frac{\sigma_{n-1}}{2} \cdot w(K)=\frac{\sigma_{n-1}}{2} \Theta$. Obviously we have also $\diam(K)= \Theta$. Hence, in our case, \eqref{eq.main.new1}
has the following form:
$$
\frac{\sigma_{n-1}}{2} \Theta=  \frac{\pi^{n/2}}{\Gamma \left(\frac{n}{2}\right)} \cdot \Theta
\leq \frac{\pi^{\frac{n-1}{2}}}{2 \Gamma \left(\frac{n+1}{2}\right)}\cdot \length (\gamma) +
\frac{\pi^{\frac{n}{2}-1}}{\Gamma \left(\frac{n}{2}\right)} \cdot \Theta
$$
or
$$
\frac{\pi-1}{\Gamma \left(\frac{n}{2}\right)} \cdot \Theta
\leq \frac{\sqrt{\pi}}{2 \Gamma \left(\frac{n+1}{2}\right)}\cdot \length (\gamma).
$$
Therefore, we get the following

\begin{theorem}\label{the3}
For a given convex body $K \subset \mathbb{R}^n$ of constant width $\Theta$
and for any  curve $\gamma \subset \mathbb{R}^n$ with the property $K\subset \co(\gamma)$, where $n\geq 2$, the inequality
\begin{equation}\label{eq.main.new2}
\length (\gamma) \geq \frac{2(\pi-1)\Gamma \left(\frac{n+1}{2}\right)}{\sqrt{\pi}\,\Gamma \left(\frac{n}{2}\right)}\cdot \Theta
\end{equation}
holds.
\end{theorem}

\begin{remark}\label{rem2}
The inequality \eqref{eq.main.new2} has the form $\length (\gamma) \geq (\pi-1)\cdot \Theta$ for  $n=2$ and
the form $\length (\gamma) \geq 4 \left(1-\frac{1}{\pi}\right)\cdot \Theta$ for $n=3$, but \eqref{eq.main.new2} does not hold in general for $n=1$.
\end{remark}

It should be noted that in the inequality  \eqref{eq.main.new2} equality does not hold if $n \geq 2$ and $\diam(K) >0$ (i.~e. $K$ is not a one-point set),
see Remark \ref{rem4}.
\smallskip

Now, we will apply Gautschi's inequality (see \cite{Gau1959})
$x^{1-s}\leq \frac{\Gamma(x+1)}{\Gamma(x+s)} \leq (x+1)^{1-s}$, where $0<s<1$ and $x>0$, to the inequality
\eqref{eq.main.new2}.
Taking $x=(n-1)/2$ and $s=1/2$, we get
\begin{equation}\label{eq.gamma1}
\sqrt{\frac{n-1}{2}}\leq \frac{\Gamma \left(\frac{n+1}{2}\right)}{\Gamma \left(\frac{n}{2}\right)} \leq\sqrt{\frac{n+1}{2}} \quad \mbox{ for all }\quad n \geq 1.
\end{equation}
This implies the following corollary from Theorem \ref{the3}:

\begin{corollary}\label{cor1}
For a given convex body $K \subset \mathbb{R}^n$ of constant width $\Theta$
and for any  curve $\gamma \subset \mathbb{R}^n$ with the property $K\subset \co(\gamma)$, the inequality
\begin{equation}\label{eq.main.new3}
\length (\gamma) \geq 2(\pi-1) \cdot \sqrt{\frac{n-1}{2\pi}}\cdot \Theta.
\end{equation}
holds.
\end{corollary}

Note also, that \eqref{eq.gamma1} implies the asymptotic
$\frac{\Gamma \left(\frac{n+1}{2}\right)}{\sqrt{n}\, \Gamma \left(\frac{n}{2}\right)} \to \frac{1}{\sqrt{2}}$ as $n \to \infty$.
\medskip

The following problem is naturally induced by Theorem \ref{the3}.

\begin{problem}\label{pr2}
For a given $n$, find the maximal constant $C(n)$ such that
$$
\length (\gamma) \geq C(n) \cdot \Theta
$$
for any convex body $K \subset \mathbb{R}^n$ of constant width $\Theta$ and for any  curve $\gamma \subset \mathbb{R}^n$ with the property $K\subset \co(\gamma)$.
\end{problem}

\section{Proof of Theorem \ref{the2} and additional remarks}\label{sect.2}

Let us recall some useful information of the Grasmann manifold $G_{n,k}$.
It could be considered as a set of $k$-dimensional unoriented linear subspaces in $\mathbb{R}^n$ with a natural manifold structure.
Note that the group $O(n)$ of all orthogonal transformation of $\mathbb{R}^n$ act transitively on $G_{n,k}$, hence,
$G_{n,k}$ is a compact homogeneous space (in fact, we have $G_{n,k}=O(n)/O(k)\times O(n-k)$).
The space $G_{n,k}$ admits $O(n)$-invariant measure, which is unique, up to a positive multiple.
For our goals, it is useful to consider the measure $\mu_{n,k}$ on $G_{n,k}$, that is introduced in \cite[II.12.4]{San2004} (although we use different notations here).
It should be noted that
$$
\int\limits_{G_{n,k}} 1\,d \mu_{n,k}=\mu_{n,k}(G_{n,k})=
\frac{\sigma_{n-1} \sigma_{n-2} \cdots \sigma_{n-k}}{\sigma_{k-1} \sigma_{k-2} \cdots \sigma_{1}\sigma_{0}}=: C_{n,k}\,.
$$
In particular, the measure $\overline{\mu_{n,k}}:=\left(C_{n,k} \right)^{-1} \cdot \mu_{n,k}$ is such that $\overline{\mu_{n,k}}(G_{n,k})=1$ and is
{\it the Haar measure} on $G_{n,k}$.

Let us fix a convex body $A$ in $\mathbb{R}^n$ and  an integer $k$, $1\leq k \leq n-1$. Now, let us consider some $P \in G_{n,k}$ and the orthogonal projection
$A'$ of $A$ to $L$. We will use notation $W'_i$ for the quermassintegrals in $P$.
If $j$ is such that $0 \leq j \leq k \leq n-1$, then we have the following equality (see \cite[(13.31)]{San2004}):

\begin{equation}\label{eq.projmes1}
\int\limits_{G_{n,k}} W'_j(A')\,d \mu_{n,k}=
\frac{n \sigma_{n-2} \cdots \sigma_{n-k}}{k \sigma_{k-2} \cdots \sigma_{0}}\cdot W_{n-k+j} (A)=
\frac{n}{k} \frac{\sigma_{k-1}}{\sigma_{n-1}}\cdot \mu_{n,k}(G_{n,k})\cdot W_{n-k+j} (A).
\end{equation}

This equality could be rewritten for the Haar measure $\overline{\mu_{n,k}}$ as follows:

\begin{equation}\label{eq.projmes2}
\int\limits_{G_{n,k}} W'_j(A')\,d\, \overline{\mu_{n,k}}=
\frac{n}{k} \frac{\sigma_{k-1}}{\sigma_{n-1}}\cdot W_{n-k+j} (A).
\end{equation}

For our goal, the most important case is when $j=k-1$:

\begin{equation}\label{eq.projmes3}
\int\limits_{G_{n,k}} W'_{k-1}(A')\,d\, \overline{\mu_{n,k}}=
\frac{n}{k} \frac{\sigma_{k-1}}{\sigma_{n-1}}\cdot W_{n-1} (A).
\end{equation}

Now, take a line segment $L$ in $\mathbb{R}^n$ as the body $A$ in \eqref{eq.projmes3}. Let us put $l=\length(L)$ and $l'=\length(L')$.
Using \eqref{eq.steiner}, we easily get
$$
\vol(L_{\varepsilon})=\omega_{n-1}\cdot l \cdot \varepsilon^{n-1}+ \omega_{n} \cdot \varepsilon^{n}\quad \mbox{ and  }\quad
\vol(L'_{\varepsilon})=\omega_{k-1}\cdot l' \cdot \varepsilon^{k-1}+ \omega_{k} \cdot \varepsilon^{k}.
$$
Consequently, $W_{n-1}(L)=\frac{\omega_{n-1}}{n} l=\frac{\sigma_{n-2}}{(n-1)n} l$ and $W'_{k-1}(L')=\frac{\omega_{k-1}}{k} l'=\frac{\sigma_{k-2}}{(k-1)k} l$.
Therefore, \eqref{eq.projmes3} implies

\begin{equation}\label{eq.projmes4}
\int\limits_{G_{n,k}} l'\,d\, \overline{\mu_{n,k}}=
\frac{(k-1)\sigma_{n-2}\sigma_{k-1}}{(n-1)\sigma_{k-2}\sigma_{n-1}}\cdot l= \frac{\sigma_{n}\sigma_{k-1}}{\sigma_{k}\sigma_{n-1}}\cdot l.
\end{equation}

The last equality is due to the relation
$$
\sigma_m=\frac{2\pi^{(m+1)/2}}{\Gamma\bigl((m+1)/2\bigr)}=\frac{4\pi^{(m+1)/2}}{(m-1)\Gamma\bigl((m-1)/2\bigr)}=\frac{2\pi}{m-1} \sigma_{m-2}
$$
for every natural $m$. Since \eqref{eq.projmes4} is fulfilled for any line segment $L$ in $\mathbb{R}^n$, then we get the equality

\begin{equation}\label{eq.projmes5}
\int\limits_{G_{n,k}} \length(\gamma')\,d\, \overline{\mu_{n,k}}= \frac{\sigma_{n}\sigma_{k-1}}{\sigma_{k}\sigma_{n-1}}\cdot \length(\gamma),
\end{equation}
where $\gamma$ is any broken line in $\mathbb{R}^n$ and $\gamma'$ is its orthogonal projection to the subspace $P\in G_{n,k}$.
Using the passage to the limit, we find that the formula \eqref{eq.projmes5} is valid for any rectifiable curve $\gamma$ in $\mathbb{R}^n$.
This result is also well known, see e.~g. \cite[Theorem 4.8.1]{AlResh1989}
\bigskip

It should be noted that the assertion of Theorem \ref{the2} is trivial if the curve $\gamma$ is not rectifiable (in this case $\length (\gamma)= \infty$).
Now we consider any convex body $K$ in $\mathbb{R}^n$ and any rectifiable curve $\gamma$ in $\mathbb{R}^n$ such that $K \subset \co(\gamma)$.
For a two-dimensional linear subspace $P$  of $\mathbb{R}^n$, we denote by $K'$ and $\gamma'$ the orthogonal projections of $K$ and $\gamma$ to $P\in G_{n,2}$
respectively. Then, as easy to see, we get $K' \subset \co(\gamma')\subset P$ and we can apply Theorem \ref{the1} in this situation.
We get the following inequality:
$$
\per (K') \leq \length (\gamma') + \diam(K')
$$

Now, let us integrate this inequality by $G_{n,k}$ with respect to the Haar measure $\overline{\mu_{n,k}}$, using the equality $\per (K')=2W'_1(K')$:

$$
2 \int\limits_{G_{n,k}} W'_1(K')\,d\, \overline{\mu_{n,k}}\leq \int\limits_{G_{n,k}} \length(\gamma')\,d\, \overline{\mu_{n,k}}+
\int\limits_{G_{n,k}}\diam(K') \,d\,  \overline{\mu_{n,k}}.
$$

Due to \eqref{eq.projmes3} and \eqref{eq.projmes5} for $k=2$, it is equivalent to the inequality

$$
\frac{\sigma_{1}}{\sigma_{n-1}}\cdot N(K)=n \frac{\sigma_{1}}{\sigma_{n-1}}\cdot W_{n-1} (K)
\leq
\frac{\sigma_{n}\sigma_1}{\sigma_{2}\sigma_{n-1}}\cdot \length(\gamma)+
\int\limits_{G_{n,k}}\diam(K') \,d\,  \overline{\mu_{n,k}}.
$$
Since $\sigma_1=2\pi$ and $\sigma_2=4\pi$, we get

\begin{equation}\label{eq.projmes6}
N(K)
\leq
\frac{\sigma_{n}}{4 \pi}\cdot \length(\gamma)+
\frac{\sigma_{n-1}}{2\pi} \int\limits_{G_{n,k}}\diam(K') \,d\,  \overline{\mu_{n,k}}.
\end{equation}

Obviously, we have $\diam(K') \leq \diam(K)$ for every $P\in G_{n,2}$, hence, we get

\begin{equation}\label{eq.projmes7}
N(K)
\leq
\frac{\sigma_{n}}{4 \pi}\cdot \length(\gamma)+
\frac{\sigma_{n-1}}{2\pi}  \diam(K).
\end{equation}

Since $\sigma_{n}=\frac{2\pi^{(n+1)/2}}{\Gamma((n+1)/2)}$ and $\sigma_{n-1}=\frac{2\pi^{n/2}}{\Gamma(n/2)}$, we get \eqref{eq.main.new1}.
Hence, we have proved Theorem~\ref{the2}.
\smallskip

\begin{remark}\label{rem3}
It should be noted that in the inequality  \eqref{eq.main.new1} equality does not hold $n \geq 3$ and $\diam(K) >0$ (i.~e. $K$ is not a one-point set).
Let suppose the contrary, i.~e. there are a convex body $K$ and a curve $\gamma$ in $\mathbb{R}^n$ for some $n \geq 3$ such that
$K \subset \co(\gamma)$ and $N(K) = \frac{\pi^{\frac{n-1}{2}}}{2 \Gamma \left(\frac{n+1}{2}\right)}\cdot \length (\gamma) +
\frac{\pi^{\frac{n}{2}-1}}{\Gamma \left(\frac{n}{2}\right)} \cdot  \diam(K)$. Now, as we can see in the above proof,
the equalities $\per (K') = \length (\gamma') + \diam(K')$  and $\diam(K')=\diam(K)$ hold, where $K'$ and $\gamma'$ are respectively
the orthogonal projections of $K$ and $\gamma$ to any two-dimensional linear subspace $P\in G_{n,2}$.
On the other hand, from Theorem \ref{the1} we get that
$\gamma'$ is a convex curve, $\bd (K')=\gamma' \cup [A',B']$, and $\diam(K')=d(A',B')$, where $A'$ and $B'$ are the endpoints of $\gamma'$.
But the latter is impossible. Indeed, $A'$ and $B'$ are the orthogonal projections of the endpoints of the curve $\gamma$
(recall that we consider the orthogonal projection from $\mathbb{R}^n$ onto $P$), hence there are two-dimensional subspaces
$P\in G_{n,2}$ with the property $B'=A'$. For any such subspace we have $0=d(A',B')=\diam(K') < \diam(K)$, that is impossible by our assumptions.
\end{remark}

\begin{remark}\label{rem4}
The arguments in Remark \ref{rem3} show that in the inequality  \eqref{eq.main.new2} equality does not hold if $n \geq 3$ and $\diam(K) >0$ (i.~e. $K$ is not a one-point set).
Moreover, the same we have also for $n=2$. Indeed, let us suppose that there are a planar convex figure of constant width $\Theta$ and a planar curve $\gamma$ such that
$K \subset \co(\gamma)$ and $\length (\gamma) = (\pi-1)\cdot \Theta$. It is clear that $\diam(K)=\Theta$ and $\per(K)=N(K)=\pi \Theta$ (see \eqref{eq.norm}).
On the other hand, from Theorem \ref{the1} we get that
$\gamma$ is a convex curve, $\bd (K)=\gamma \cup [A,B]$, and $\diam(K)=d(A,B)$, where $A$ and $B$ are the endpoints of $\gamma$.
We know that the boundary of any planar figure of constant width contains no line segment (i.~e. a figure of constant width is strictly convex), see e.~g.
Theorem 3.1.1 and the discussion after its statement in \cite{MaMoOl}.
Hence, $B=A$ and $\diam(K)=d(A,B)=0$.
\end{remark}


\vspace{10mm}

\begin{thebibliography}{99}

\bibitem{AlResh1989}
A.D.~Alexandrov, Yu.G.~Reshetnyak, {\sl General theory of irregular curves},
Mathematics and its Applications (Soviet Series), 29.
Kluwer Academic Publishers Group, Dordrecht, 1989. Translated from the Russian by L.Ya.~Yuzina, {\bf Zbl.}0691.53002, {\bf MR}1117220.

\bibitem{AkVys2017}
A.~Akopyan, V.~Vysotsky,
{\sl On the lengths of curves passing through boundary points of a planar convex shape,}
Am. Math. Mon., 124(7) (2017), 588--596, {\bf Zbl.}1391.52003,  {\bf MR}3681589.

\bibitem{BoFe1987}
T.~Bonnesen, W.~Fenchel,
{\sl Theory of convex bodies}, BCS Associates, Moscow, ID, 1987. Translated
from the German and edited by L. Boron, C. Christenson and B. Smith., {\bf Zbl.}0628.52001, {\bf MR}0920366.

\bibitem{BBI2001}
D. Burago, Yu. Burago, and S. Ivanov.
{\sl A course in metric geometry}, Graduate Studies in Mathematics, 33.
American Mathematical Society, Providence, RI, 2001, {\bf Zbl.}0981.51016, {\bf MR}1835418.

\bibitem{Gau1959}
W.~Gautschi,
{\sl Some elementary inequalities relating to the gamma and incomplete gamma function,}
J. Math. Phys., 38 (1959), 77--81, {\bf Zbl.}0094.04104,  {\bf MR}0103289.


\bibitem{Ghomi2018}
M.~Ghomi,
{\sl The length, width, and inradius of space curves,}
Geom. Dedicata, 196 (2018), 123--143, {\bf Zbl.}1403.53005,  {\bf MR}3853631.


\bibitem{Had1957}
H.~Hadwiger,
{\sl Vorlesungen \"uber Inhalt, Oberfl\"ache und Isoperimetrie,}
Die Grundlehren der Mathematischen Wissenschaften, 93. Springer-Verlag, Berlin--G\"{o}ttingen--Heidelberg, 1957, {\bf Zbl.}0078.35703,  {\bf MR}0102775.

\bibitem{MaMoOl}
H.~Martini, L.~Montejano, D.~Oliveros,
{\sl Bodies of constant width. An introduction to convex geometry with applications,}
Birkh\"auser/Springer, Cham, 2019, {\bf MR}3930585, {\bf Zbl.}06999635.


\bibitem{NikNik2021}
Yu.G.~Nikonorov, Yu.V.~Nikonorova,
{\sl One property of a planar curve whose convex hull covers a given convex figure,}
Elemente der Mathematik, 2021, DOI: 10.4171/EM/458, see also \textit{arXiv}:2007.00612.


\bibitem{San2004}
L.A.~Santal\'{o},
{\sl Integral geometry and geometric probability,}
Second edition, with a foreword by M.~Kac, Cambridge Mathematical Library. Cambridge University Press, Cambridge, 2004,
{\bf Zbl.}1116.53050,  {\bf MR}2162874.

\bibitem{Zal1996}
V.A.~Zalgaller, {\sl Extremal problems on the convex hull of a space curve (Russian)},
Algebra i Analiz, 8(3) (1996),  1--13, English translation: St. Petersburg Math. J., 8(3) (1997), 369--379,  {\bf Zbl.}0877.52005,  {\bf MR}1402285.


\end{thebibliography}
\end{document}